\documentclass{siamltex}
\usepackage{graphicx}
\usepackage{amsfonts,latexsym}
\usepackage{amsmath}
\newcommand{\R}{{\mathbb R}}
\newcommand{\C}{{\mathbb C}}

\newcommand{\bx}{\mbox{\boldmath{$x$}}}
\newcommand{\bb}{\mbox{\boldmath{$b$}}}

\newcommand{\br}{\mbox{\boldmath{$r$}}}

\newcommand{\be}{\mbox{\boldmath{$e$}}}

\newcommand{\by}{\mbox{\boldmath{$y$}}}

\newcommand{\bbeta}{\mbox{\boldmath{$\beta$}}}
\newcommand{\bxi}{\mbox{\boldmath{$\xi$}}}
\newcommand{\bdxi}{\mbox{\boldmath{$\Delta\xi$}}}
\newcommand{\bdbeta}{\mbox{\boldmath{$\Delta\beta$}}}

\newcommand{\sbx}{\mbox{\boldmath{${\scriptstyle x}$}}}

\newcommand{\sbr}{\mbox{\boldmath{${\scriptstyle r}$}}}

\newcommand{\bzero}{\mbox{\boldmath{$0$}}}
\newcommand{\bl}{\begin{list}{ \ }{
\leftmargin=.325in}}
\newcommand{\el}{\end{list}}

\begin{document}
\date{}
\title{Circulant preconditioners for \\ discrete ill-posed Toeplitz systems}
\author{L. Dykes\thanks{Department of Mathematical Sciences, Kent State University, Kent, 
OH 44242, USA, and University School, Hunting Valley, OH 44022, USA.  E-mail: 
{\tt ldykes@math.kent.edu}.}\and
S. Noschese\thanks{Dipartimento di Matematica ``Guido Castelnuovo'', SAPIENZA Universit\`a
di Roma, P.le A. Moro, 2, I-00185 Roma, Italy. E-mail: {\tt noschese@mat.uniroma1.it}.
Research supported by a grant from SAPIENZA Universit\`a di Roma.}\and
L. Reichel\thanks{Department of Mathematical Sciences, Kent State University, Kent, OH 
44242, USA. E-mail: {\tt reichel@math.kent.edu}.}
}
\maketitle
\begin{center}
{\em Dedicated to Ken Hayami on the occasion of his 60th birthday.}
\end{center}

\begin{abstract}
Circulant preconditioners are commonly used to accelerate the rate of convergence of
iterative methods when solving linear systems of equations with a Toeplitz matrix. Block
extensions that can be applied when the system has a block Toeplitz matrix with Toeplitz 
blocks also have been developed. This paper is concerned with preconditioning of linear
systems of equations with a symmetric block Toeplitz matrix with symmetric Toeplitz blocks 
that stem from the discretization of a linear ill-posed problem. The right-hand side of 
the linear systems represents available data and is assumed to be contaminated by error. 
These kinds of linear systems arise, e.g., in image deblurring problems. It is important 
that the preconditioner does not affect the invariant subspace associated with the 
smallest eigenvalues of the block Toeplitz matrix to avoid severe propagation of the error
in the right-hand side. A perturbation result indicates how the dimension of the subspace 
associated with the smallest eigenvalues should be chosen and allows the determination of
a suitable preconditioner when an estimate of the error in the right-hand side is 
available. This estimate also is used to decide how many iterations to carry out by a 
minimum residual iterative method. Applications to image restoration are presented.
\end{abstract}



\section{Introduction}\label{sec1}
Linear systems of equations with a matrix with a Toeplitz-type structure arise in many 
applications, such as in signal and image processing. Consider the computation of an 
approximate solution of the linear system of equations
\begin{equation}\label{linsys}
T{\bx}={\bb}, \qquad T\in{\R}^{n_1n_2\times n_1n_2},\qquad {\bx},{\bb}\in{\R}^{n_1n_2},
\end{equation}
where $T$ is a symmetric BTTB matrix, i.e., $T$ is a symmetric block Toeplitz matrix with each block 
being an $n_1 \times n_1$ symmetric Toeplitz matrix. The eigenvalues of $T$ are assumed to
decay smoothly to zero in magnitude without a significant gap. In particular, $T$ may be 
singular. Linear systems of equations (\ref{linsys}) with a matrix of this kind arise, for example, 
from the discretization of a linear ill-posed problem, such as a Fredholm integral 
equation of the first kind in two space-dimensions with a displacement kernel. 

The right-hand side $\bb$ of (\ref{linsys}) is assumed to be contaminated by an (unknown) 
error $\be$. We will refer to this error as ``noise''. It may stem from measurement or 
discretization errors. Let $\widehat{\bb}$ denote the (unknown) error-free vector 
associated with ${\bb}$, i.e., 
\begin{equation}\label{rhs}
{\bb}=\widehat{\bb}+{\be}.
\end{equation}
The (unknown) linear system of equations with error-free right-hand side,
\begin{equation}\label{linsysnf}
T{\bx}=\widehat{\bb},
\end{equation}
is assumed to be consistent; however, we do not require the available system 
(\ref{linsys}) to be consistent.

We will assume that a fairly sharp bound for the norm of $\be$ is known, i.e.,
\begin{equation}\label{errbd}
\|\be\|\leq\varepsilon.
\end{equation}
Here and throughout this paper $\|\cdot\|$ denotes the Euclidean vector norm or the 
spectral matrix norm. This bound will help us determine a suitable number of iterations 
to carry out with a minimal residual iterative method and to construct a preconditioner 
for the solution of \eqref{linsys}.

Let $T^\dag$ denote the Moore--Penrose pseudoinverse of $T$. We are interested in 
computing an approximation of the solution $\widehat{\bx}=T^\dag\widehat{\bb}$ of minimal 
Euclidean norm of the unavailable error-free linear system (\ref{linsysnf}). Note that the 
solution of (\ref{linsys}),
\[
\bx=T^\dag{\bb}=T^\dag(\widehat{\bb}+{\be})=\widehat{\bx}+T^\dag{\be},
\]
typically is dominated by the propagated error $T^\dag{\be}$ and, therefore, is useless.
Therefore all solution methods for \eqref{linsys} seek to determine a suitable approximate
solution that is not severely contaminated by propagated error. The computed approximate 
solution is the exact solution of an appropriately chosen nearby problem, whose solution
is less sensitive to the error $\be$ in $\bb$ than the solution of \eqref{linsys}. The 
replacement of the given problem 
\eqref{linsys} by a nearby problem is commonly referred to as regularization. Among the
most popular regularization methods is Tikhonov regularization, which replaces 
\eqref{linsys} by a penalized least-squares problem, and truncated iteration, which is
based on solving \eqref{linsys} by an iterative method and terminating the iterations
suitably early, see, e.g., \cite{EHN,Ha0,Ki,RR} for discussions on these regularization 
methods. In this paper we regularize by truncated iteration, and by choosing a suitable 
preconditioner. 

The evaluation of matrix-vector products with a BTTB matrix of order $n_1n_2$ can be 
carried out in only ${\mathcal O}(n_1n_2\log_2(n_1n_2))$ arithmetic floating point 
operations (flops) by using the fast Fourier transform (FFT); see, e.g., \cite{CJ,Ng}. 
This makes it attractive to solve \eqref{linsys} by an iterative method. We will use a
preconditioner to increase the rate of convergence of the iterative method. BTTB matrices 
are commonly preconditioned by block circulant with circulant block (BCCB) matrices; see
\cite{CJ,CO,Ng,vdMRS1,vdMRS2} for discussions, illustrations, and further references. The 
use of BCCB preconditioners for BTTB matrices is attractive both due to the spectral 
properties of the preconditioned matrix and because of the possibility to evaluate a 
matrix-vector product with a preconditioned matrix of order $n_1n_2$ in only 
${\mathcal O}(n_1n_2\log_2(n_1n_2))$ flops with the aid of the FFT. A MATLAB software
package for fast matrix-vector product evaluation is provided by Redivo--Zaglia and 
Rodriguez \cite{RZR}.

When preconditioning a BTTB matrix $T$ that stems from the discretization of a linear 
ill-posed problem, it is desirable that an invariant subspace associated with the
eigenvalues of $T$ of smallest magnitude is not affected by preconditioning to avoid 
severe propagation of the error $\be$ in the right-hand side $\bb$ of (\ref{linsys}) into
the computed iterates. This is due to the fact that the eigenvectors associated with these
eigenvalues are highly oscillatory (have many sign changes) and model noise rather than
the desired solution $\widehat{\bx}$. Typically, we do not want these eigenvectors to be
part of our computed approximation of $\widehat{\bx}$. A nice introduction to BCCB
preconditioners for the solution of discretized linear ill-posed problems with a 
Toeplitz-type matrix is presented by Hanke et al. \cite{HNP}. 

The number of iterations have to be few enough to avoid severe propagation of the error
$\be$ in $\bb$ into the computed approximation of $\widehat{\bx}$. The availability of 
the bound \eqref{errbd} and the consistency of \eqref{linsysnf} allow us to apply the
discrepancy principle to determine a suitable number of iterations as well as to 
define the dimension of the invariant subspace that should not be affected by the
preconditioner. Roughly, the larger the error in $\bb$, the larger should the dimension of
the subspace that is not (or only minimally) affected by preconditioning be chosen. 

Various other approaches to define BCCB preconditioners for the iterative solution of
discretized linear ill-posed problems \eqref{linsys} with an error-contaminated 
right-hand side are described in the literature. For instance, Hanke and Nagy \cite{HN}
apply the L-curve criterion to determine a subspace that should not be affected by the
 preconditioner. The L-curve criterion implicitly estimates the norm of the error in 
$\bb$. This criterion is able to estimate the norm of the error fairly accurately in 
some situations, but it is not a reliable error estimator; see \cite{Ki,RR} for 
discussions and illustrations. We therefore are interested in developing an approach 
for constructing BCCB preconditioners that is not based on the L-curve criterion. Hanke
et al. \cite{HNP} apply a discrete Picard condition to determine the dimension of the 
subspace that the BCCB preconditioner should leave invariant. This approach typically works
quite well in an interactive computing environment that allows the determination of whether
the discrete Picard condition holds by visual inspection, however, it is not straightforward
to automatize. Di Benedetto et al. \cite{DBES} propose the application of a so-called 
superoptimal BCCB preconditioner and do not explicitly choose the dimension of the subspace
that should be unaffected by the preconditioner. This approach works well for some image
restoration problems, but not for others; see the discussion in \cite{DBES}. 

Preconditioning is most useful when the error $\be$ in $\bb$ is of small relative
norm, because then many steps of an iterative method may be required to determine an
accurate approximation of $\widehat{\bx}$. When the error $\be$ is large, only 
few steps can be carried out before the propagated error destroys the computed solution.
Preconditioning then does not reduce the computational effort by much. 

In this paper, we will use the bound (\ref{errbd}) to determine both the BCCB 
preconditioner and the number of iterations to be carried out. A perturbation bound guides
our choice of preconditioner. This is described in Section \ref{sec2}. A few computed 
examples are presented in Section \ref{sec3} and concluding remarks can be found in 
Section \ref{sec4}.

\section{Preconditioned iterative regularization}\label{sec2}
We discuss the construction of the preconditioner, the stopping criterion for the 
iterative method, and outline the minimal residual iterative methods used.

\subsection{The BCCB preconditioner}\label{sub2.1}
Let for the moment $T\in\R^{n\times n}$ be a symmetric positive definite Toeplitz matrix
and let $C\in\R^{n\times n}$ be the closest circulant matrix to $T$ in the Frobenius norm.
T. Chan \cite{TC} proposed the use of $C$ as a preconditioner for $T$; see also 
\cite{CJ,Ng}. The eigenvalues of $C$ are given by the discrete Fourier transform of the 
first column of $C$; their computation with the FFT requires only 
${\mathcal O}(n\log_2(n))$ flops. Since we would like the preconditioner not to affect the 
invariant subspace of $T$ associated with the smallest eigenvalues, we set the $n-p$ 
eigenvalues of smallest magnitude of $C$ to unity for some suitable $0\leq p\leq n$, analogously
as in \cite{HN,HNP}. We refer to this preconditioner as $C_{p}$. Subsection \ref{sub2.4} 
describes how to determine $p$ using the error bound \eqref{errbd}. 

The linear systems (\ref{linsys}) of interest to us have a BTTB matrix $T$, i.e., $T$ is
the Kronecker product of two Toeplitz matrices 
\begin{equation}\label{TT}
T=T_1\otimes T_2,\qquad T_1\in\R^{n_1\times n_1},\qquad T_2\in\R^{n_2\times n_2}.
\end{equation}
We will determine a preconditioner that is the Kronecker product of two circulant matrices
\begin{equation}\label{CC}
C=C_{p_1}\otimes C_{p_2},\qquad C_{p_1}\in\R^{n_1\times n_1},\qquad 
C_{p_2}\in\R^{n_2\times n_2},
\end{equation}
where $C_{p_j}$ is defined by first determining the closest circulant $C_j$ to $T_j$ in the
Frobenius norm and then setting its $n_j-p_j$ eigenvalues of smallest magnitude to one for
$j=1,2$. In this way our preconditioner $C$ does not affect the invariant subspace of $T$ associated with the eigenvalues of smallest magnitude. The eigenvectors of this subspace are highly oscillatory
and primarily model noise and not the desired solution $\widehat{\bx}$. This construction of $C$ requires only ${\mathcal O}(n_1n_2(\log_2(n_1)+\log_2(n_2))$ flops. We refer to the BCCB 
preconditioner so determined as $C_{p_1,p_2}$. The determination of this preconditioner is
somewhat faster than of the BCCB preconditioner described in \cite{HNP}, because the latter
requires that all its $n_1n_2$ eigenvalues be formed.

\subsection{Stopping criterion}\label{sub2.3}
Once the BCCB preconditioner $C_{p_1,p_2}$ has been defined, we compute an approximate 
solution ${\by}_k$ of the preconditioned linear system of equations,
\begin{equation}
\label{prec_sys}
TC_{p_1,p_2}^{-1}{\by}={\bb},
\end{equation}
using one of the minimal residual iterative methods described in Subsection \ref{sub2.6}. 
As for the initial approximation of $\widehat{\bx}$, define the BCCB matrix 
$\widetilde{C}_{p_1,p_2}=\widetilde{C}_{p_1}\otimes \widetilde{C}_{p_2}$, where 
$\widetilde{C}_{p_j}$ is obtained from $C_j$ by setting the $n_j-p_j$ eigenvalues of 
smallest magnitude to zero, for $j=1,2$, and define 
\begin{equation}
\label{init_guess} 
{\bx}_0=\widetilde{C}_{p_1,p_2}^{\dag}{\bb},\qquad {\br}_0={\bb}-T{\bx}_0. 
\end{equation}
The initial iterate then is chosen to be ${\by}_0= C_{p_1,p_2}{\bx_0}$. 

Let ${\by}_1,{\by_2},\ldots~$ denote the computed iterates. The number of iterations to be
carried out is determined with the aid of the discrepancy principle. This stopping 
criterion prescribes that the iterations be terminated as soon as an iterate $\by_k$ 
satisfies
\begin{equation}\label{discrp1}
\|TC_{p_1,p_2}^{-1}\by_k-\br_0\|\leq\gamma\varepsilon,
\end{equation}
where $\gamma \geq 1$ is a user-specified parameter independent of $\varepsilon$. 
Typically, $\gamma$ is chosen close to unity when $\varepsilon$ is known to be a fairly 
sharp upper bound for $\|\be\|$; cf. \eqref{errbd}. We obtain the approximation 
\[
{\bx}_k={\bx_0}+C_{p_1,p_2}^{-1}{\by}_k
\]
of the desired vector $\widehat{\bx}$.

\subsection{Construction of the preconditioner}\label{sub2.4}
Let for the moment $T\in\R^{n\times n}$ be a Toeplitz matrix and let $C\in\R^{n\times n}$
be the closest circulant in the Frobenius norm. Order the eigenvalues of $C$ according to
\[
|\lambda_1|\geq |\lambda_2|\geq \dots \geq|\lambda_n|\geq 0.
\]
Assume that a bound (\ref{errbd}) is known. Let $p$ be the number of eigenvalues of 
largest magnitude of $C$ that are not set to unity. We choose 
$p=\lfloor\frac{3}{4}q\rfloor$, where $q$ is the solution of the minimization problem
\begin{equation}\label{choice_p}
\min_{1\leq q<n}\frac{1}{|\lambda_q|}
\left(\frac{|\lambda_{q+1}|}{|\lambda_{1}|}+\eta\right).
\end{equation}
Here $\lfloor\alpha\rfloor$ denotes the largest integer smaller than or equal to 
$\alpha\geq 0$ and $\eta=\varepsilon/\|{\bb}\|$. This choice of $p$ is suggested by 
the following perturbation result.

\begin{proposition}\label{prop1}  
Given a rank-$q$ matrix $A\in{\C}^{n\times n}$, $q\leq n$, a vector $\bbeta\ne\bzero$ 
in the range of $A$, and ${\bxi}$ such that ${\bxi} = A^{\dag}{\bbeta}$. Let $\Delta A$,
$\bdbeta$, and ${\bdxi}$ satisfy 
\[
(A+\Delta A)({\bxi}+{\bdxi})= {\bbeta}+{\bdbeta}.
\]
Then
\begin{equation}
\label{ineq}
\frac{\|{\bdxi}\|}{\|{\bxi}\|}
\leq  \nu(A,A+\Delta A)\, \kappa(A)
\left( \frac{\|\Delta A\|}{\|A\|}+\frac{\|{\bdbeta}\|}{\|{\bbeta}\|}\right),
\end{equation}
where $\kappa(M)=\|M\|\|M^{\dag}\|$ is the condition number of the matrix $M$, and 
$\nu(M,N)$ denotes the ratio between the smallest singular values of the matrices $M$ and
$N$.
\end{proposition}

\proof
From $A{\bxi} ={\bbeta}$ and $(A+\Delta A){\bdxi} = -\Delta A {\bxi} + {\bdbeta}$, one has  
\[
{\bdxi}=(A+\Delta A)^{\dag}\left( -\Delta A{\bxi}+{\bdbeta}\right).
\]
Taking norms on both sides, one gets
\[
\|{\bdxi}\| \leq \|(A+\Delta A)^{\dag}\| \left( \|\Delta A\|\|{\bxi}\| +\|{\bdbeta}\|\right).
\]
Hence,
\[
\|{\bdxi}\| \leq \nu(A,A+\Delta A)\, \kappa(A)
\left(\frac{\|\Delta A\|\|{\bxi}\| }{\|A\|}+\frac{\|\bdbeta\|}{\|A\|}\right).
\]
Finally, dividing by $\|{\bxi}\|$ and exploiting the inequality 
$\|{\bbeta}\|\leq\|A\|\|\bxi\|$ yield (\ref{ineq}). 
\endproof

Let $\widetilde{C}_{q}$ be the circulant obtained by setting the $n-q$ eigenvalues of 
smallest magnitude of $C$ to zero. In this context, one replaces $A$ by $\widetilde{C}_{q}$ 
and ${\bbeta}$ by ${\bb}$ in Proposition \ref{prop1}. Now, letting 
$\Delta A=T-\widetilde{C}_{q}$ and ${\bdbeta}=-\be$, we obtain by (\ref{rhs}) that 
${\bdxi}=\hat{\bx}-\widetilde{C}_{q}^{\dag}{\bb}$ satisfies the 
hypothesis of Proposition \ref{prop1}. Inequality (\ref{ineq}) reads
\begin{eqnarray*}
\frac {\|\hat{\bx}-\widetilde{C}_{q}^{\dag}{\bb}\|}{\|{\bx_0}\|}
&\leq& \nu(\widetilde{C}_{q},T)\, \kappa(\widetilde{C}_{q})
\left(\frac{\|T-\widetilde{C}_{q}\|}{\|\widetilde{C}_{q}\|}+
\frac{\|\be\|}{\|\bb\|}\right) \\
&\leq&  \nu(\widetilde{C}_{q},T)\, \frac{|\lambda_1|} {|\lambda_q|} 
\left(\frac{\|T-\widetilde{C}_{q}\|}{|\lambda_1|}+\eta\right) \sim 
\frac{|\lambda_1|} {|\lambda_q|} \left(\frac{|\lambda_{q+1}|}{|\lambda_1|}+\eta\right).
\end{eqnarray*}
In the final estimate, we assume that the eigenvalues of $T$ of smallest magnitude are  
close to $|\lambda_q|$, so that $\nu(\widetilde{C}_{q},T)\sim 1$ and 
$\|T-\widetilde{C}_{q}\|\sim\|C-\widetilde{C}_{q}\|$. This discussion suggests the choice
$p=q$, where $q$ is determined by (\ref{choice_p}). However, since we do not know whether
the eigenvalues of $T$ of smallest magnitude are close to $|\lambda_q|$, we will choose 
$p=\lfloor 3q/4\rfloor$ to secure that we do not ``over-precondition'' and thereby obtain 
a large propagated error in the computed approximation of $\widehat{\bx}$. We remark that
a ``standard'' circulant preconditioner $\widetilde{C}_n$ generally over-preconditions and
gives a large propagated error in the computed solution.

We turn to BTTB matrices and first consider the matrix $T\otimes T$. Let $C$ be the 
closest circulant to $T$ in the Frobenius norm. Then $C\otimes C$ is the closest BCCB 
matrix to $T\otimes T$ in the Frobenius norm. Let $q$ to be the solution
of the minimization problem
\begin{equation}\label{choice2_p}
\min_{1\leq q<n}\frac{1}{|\lambda_q|^2}\left(\frac{|\lambda_{q+1}|^2}{|\lambda_{1}|^2}+
\eta\right),
\end{equation}
where $\eta=\varepsilon/\|{\bb}\|$. The choice $p=q$ is suggested by Proposition 
\ref{prop1}, where one replaces $A$ by $\widetilde{C}_{p,p}$, the BCCB matrix obtained 
setting the $n-p$ eigenvalues of smallest magnitude of each circulant matrix $C$ to zero,
and replaces ${\bbeta}$ by ${\bb}$, so that the vector ${\bxi}$ in Proposition \ref{prop1}
is given by ${\bx_0}$ in (\ref{init_guess}) with $p_1=p_2=q$. In the computed examples, 
we will let $p_1=p_2=\lfloor 3q/4\rfloor$ to avoid to over-precondition.

Finally, consider BTTB matrices of the form (\ref{TT}). To determine a BCCB 
preconditioner of the kind (\ref{CC}), we sort the eigenvalues of the circulant matrices 
$C_j$, $j=1,2$, according to
\[
|\lambda_1^{(j)}|\geq |\lambda_2^{(j)}|\geq \dots \geq|\lambda_{n_j}^{(j)}|\geq 0,
\]
and let $p_j$ be the number of eigenvalues of largest magnitude of $C_j$ that are not set 
to unity. Let the index pair $\{q_1,q_2\}$ solve the minimization problem
\begin{equation}\label{choice3_p}
\min_{\substack{1\leq q_1<n_1\\1\leq q_2<n_2}}
\frac{1}{|\lambda^{(1)}_{q_1}||\lambda^{(2)}_{q_2}|}
\left(\frac{|\lambda^{(1)}_{q_1+1}||\lambda^{(2)}_{q_2+1}|}
{|\lambda^{(1)}_{1}||\lambda^{(2)}_{1}|}+\eta\right),
\end{equation}
where $\eta=\varepsilon/\|{\bb}\|$. Similarly as above, we let $p_j=\lfloor 3q_j/4\rfloor$
for $j=1,2$.

\subsection{Construction of the preconditioned system}\label{sub2.5} 
The preconditioned matrix $T C_{p_1,p_2}^{-1}$, where $T$ and $C_{p_1,p_2}$ are given by (\ref{TT}) and (\ref{CC}), respectively, is constructed as follows: 

\noindent\vskip3pt
Compute for $j=1,2$:
\begin{enumerate}
\item  
The closest circulant $C_j$ to $T_j$ in the Frobenius norm.
\item  
The FFT of the first column of the matrix $C_j$. This gives the eigenvalues 
$\lambda^{(j)}_1,\lambda^{(j)}_2,\ldots,\lambda^{(j)}_{n_j}$ of $C_j$. The eigenvectors are 
the columns of the Fourier matrix. Permute the columns of the Fourier matrix so that the
eigenvalues are ordered according to decreasing magnitude,  
\[
|\lambda^{(j)}_1|\ge|\lambda^{(j)}_2|\ge\ldots\ge|\lambda^{(j)}_{n_j}|
\]
and denote the permuted Fourier matrix by $U_j$. The columns of $U_j$ generally become 
more oscillatory with increasing column number.
\item  
The truncation index $p_j$, $1\leq p_j \leq n_j$, by using (\ref{choice3_p}). This defines
the diagonal matrices 
\begin{eqnarray*}
\Lambda_{p_j}&=&{\rm diag}[\lambda^{(j)}_1,\lambda^{(j)}_2,\ldots,\lambda^{(j)}_{p_j},1,\ldots,1],\\
 \widetilde{\Lambda}_{p_j}&=&{\rm diag}[\lambda^{(j)}_1,\lambda^{(j)}_2,\ldots,
\lambda^{(j)}_{p_j}, 0,\ldots, 0].
\end{eqnarray*}
Neither the matrices $U_j$ nor $\Lambda_{p_j}$ and $\widetilde{\Lambda}_{p_j}$ have to be
explicitly formed.
\end{enumerate}

We are now in a position to define the preconditioner and related matrices, but hasten 
to point out that these matrices do not have to be explicitly formed. Introduce
\[
\begin{array}{rclcrclcrcl}
C_{p_j}&=&U_j \Lambda_{p_j}U_j^*,&~~& C_{p_j}^{-1}&=&U_j\Lambda_{p_j}^{-1}U_j^*,&~~&
\widetilde{C}_{p_j}^{\dag}&=&U_j \widetilde{\Lambda}_{p_j}^{\dag} U_j^*,\\
C_{p_1,p_2}&=&C_{p_1}\otimes C_{p_2},&~~& C_{p_1,p_2}^{-1}&=&C_{p_1}^{-1}\otimes C_{p_2}^{-1},
&~~& \widetilde{C}_{p_1,p_2}^{\dag}&=&\widetilde{C}_{p_1}^{\dag}\otimes \widetilde{C}_{p_2}^{\dag},
\end{array}
\]
where the superscript $^*$ denotes transposition and complex conjugation. We compute the initial
approximate solution ${\bx_0}=\widetilde{C}_{p_1,p_2}^{\dag}{\bb}$ in (\ref{init_guess}) without
explicitly forming the matrix $\widetilde{C}_{p_1,p_2}^{\dag}$. Indeed, the spectral factorization 
\[
C_{p_1,p_2}=C_{p_1}\otimes C_{p_2}=(U_1\otimes
U_2)(\Lambda_{p_1}\otimes\Lambda_{p_2})(U_1\otimes
U_2)^*,
\]
can be applied to evaluate $C_{p_1,p_2}\by$ for any $\by\in{\R}^{n_1n_2}$ in 
${\mathcal O}(n_1n_2\log_2(n_1n_2))$ flops with the FFT, and  the same holds for
matrix-vector products with the matrices $C_{p_1,p_2}^{-1}$ and $\widetilde{C}_{p_1,p_2}^{\dag}$.

Krylov subspace methods for the iterative solution of \eqref{prec_sys} require matrix-vector
product evaluations with the preconditioned matrix $TC^{-1}_{p_1,p_2}$. It is well known that
these matrix-vector product evaluations can be carried out quickly with the aid of the FFT. We
outline for completeness the evaluation of matrix-vector products with the
matrix $TC_p^{-1}$ in the simplified situation when $T \in \R^{n\times n}$ is a Toeplitz matrix
and $C_p$ is a circulant. We express $T$ as a sum of a circulant $C_0$ and a skew-circulant 
$C_{\pi}$. This splitting and the spectral factorizations
\[
C_0=U \Lambda_0 U^*,\quad C_\pi=U_\pi  \Lambda_\pi U_\pi^*,
\]
where $U_\pi=\diag\{1,e^{i\pi/n},\dots,e^{(n-1)i\pi/n} \}U$,  yield that
\[
TC_p^{-1}=(C_0+C_\pi)C_p^{-1}=
(U \Lambda_0 U^* +U_\pi  \Lambda_\pi U_\pi^*)U \Lambda_p^{-1} U^*.
\]
The preconditioned linear system of equations $TC_p^{-1}\bx=\bb$ with $\bb \in \R^{n}$
can be expressed in the form
\[
(\Lambda_0+U^*U_\pi\Lambda_\pi U_\pi^*U)\Lambda_p^{-1}\by=U^* \bb,
\qquad \bx= U  \by,
\]
which is used in the computations. Each iteration requires the evaluation of 
the FFT of four $n$-vectors. The computation of these FFTs is the dominating 
computational work. We remark that the dominating computational effort to 
evaluate a matrix-vector product with the matrix $T$, which is required when solving 
the unpreconditioned system $T\bx=\bb$ by a Krylov subspace method, also is the
computation of the FFT of four $n$-vectors. Therefore, the number of iterations
required by the iterative method is the proper measure of the computational effort 
both for preconditioned and unpreconditioned linear systems of equations. The 
situation is analogous when $T$ is the tensor product of two Toeplitz matrices.
We omit the details.

\subsection{Range restricted GMRES and MINRES methods}\label{sub2.6} 
GMRES is a popular iterative method for the solution of large linear systems of equations 
$A{\bx}={\bb}$ with a square nonsingular nonsymmetric matrix $A$ that arise from the 
discretization of well-posed problems, such as Dirichlet boundary value problems for 
elliptic partial differential equations; see, e.g., Saad \cite{Sa}. The $k$th iterate 
determined by this method solves the minimization problem
\[
\min_{\sbx\in{\mathcal K}_k(A,{\sbr}_0)}\|A\bx-\br_0\|,
\]
where $\bx_0$ is an initial approximate solution, ${\br}_0={\bb}-A{\bx}_0$, and 
\[
{\mathcal K}_k(A,{\br}_0)=\mbox{span}\{\br_0,A\br_0,\ldots,A^{k-1}{\br}_0\}
\]
is a Krylov subspace. 

It has been observed that a modification of GMRES, which we refer to as the range 
restricted GMRES method (RRGMRES), often yields a more accurate approximation of the 
desired solution $\widehat{\bx}$ than (standard) GMRES when $A$ stems from the 
discretization of a linear ill-posed problem
and the right-hand side $\bb$ is contaminated by error; see \cite{DR,NRS}. The
RRGMRES method determines iterates in shifted Krylov subspaces 
${\mathcal K}_k(A,A^\ell{\br}_0)$, where $\ell\geq 1$ is a small integer. We propose
that an RRGMRES method be used for the solution of the preconditioned problem
\[
TC_{p_1,p_2}^{-1}\by=\br_0.
\]
When $T$ and $C_{p_1,p_2}$ are symmetric positive definite, RRGMRES can be simplified
to a range restricted MINRES method that only requires simultaneous storage of a few
$n$-vectors, the number of which is bounded independently of the number of iterations;
see \cite{DMR} for details.

\section{Computed examples}\label{sec3}
The calculations of this section were carried out in MATLAB with machine epsilon about 
$2.2\cdot 10^{-16}$. For all the examples, we chose  $\gamma=1$ in (\ref{discrp1}).

{\bf Example 1.} This is an image deblurring test problem from the MATLAB package 
Regularization Tools \cite{Tools}. The original  image and the symmetric BTTB  matrix 
are those defined by the MATLAB function {\sf blur.m}. We choose the dimensions 
$n_1=n_2=64$, half-bandwidth for each Toeplitz block specified by the parameter $band=10$
and width of the Gaussian point spread specified by the parameter $sigma=\sqrt5$.

\begin{table}[htb!]
\centering
\begin{tabular}{cccc}\hline
\% relative data error & $p$ & steps $k$ & $\|\bx_k-\widehat{\bx}\|/\|\widehat{\bx}\|$ \\
\hline
$0.10$ & $14$ & $18$ & $0.3404$ \\
$0.10$ & $-$ & $33$ & $0.3361$ \\ \hline
$0.05$ & $16$ & $22$ & $0.3308$ \\
$0.05$ & $-$ & $45$ & $0.3275$ \\ \hline
$0.01$ & $17$ & $42$ & $0.3094$ \\
$0.01$ & $-$ & $89$ & $0.3072$ \\ \hline
\end{tabular}
\caption{Example 1: blur. A hyphen signifies that no preconditioner is used.}
\label{tab02}
\end{table}

\begin{figure}[ht]
\centering
\begin{tabular}{cc}
\includegraphics[width=6cm]{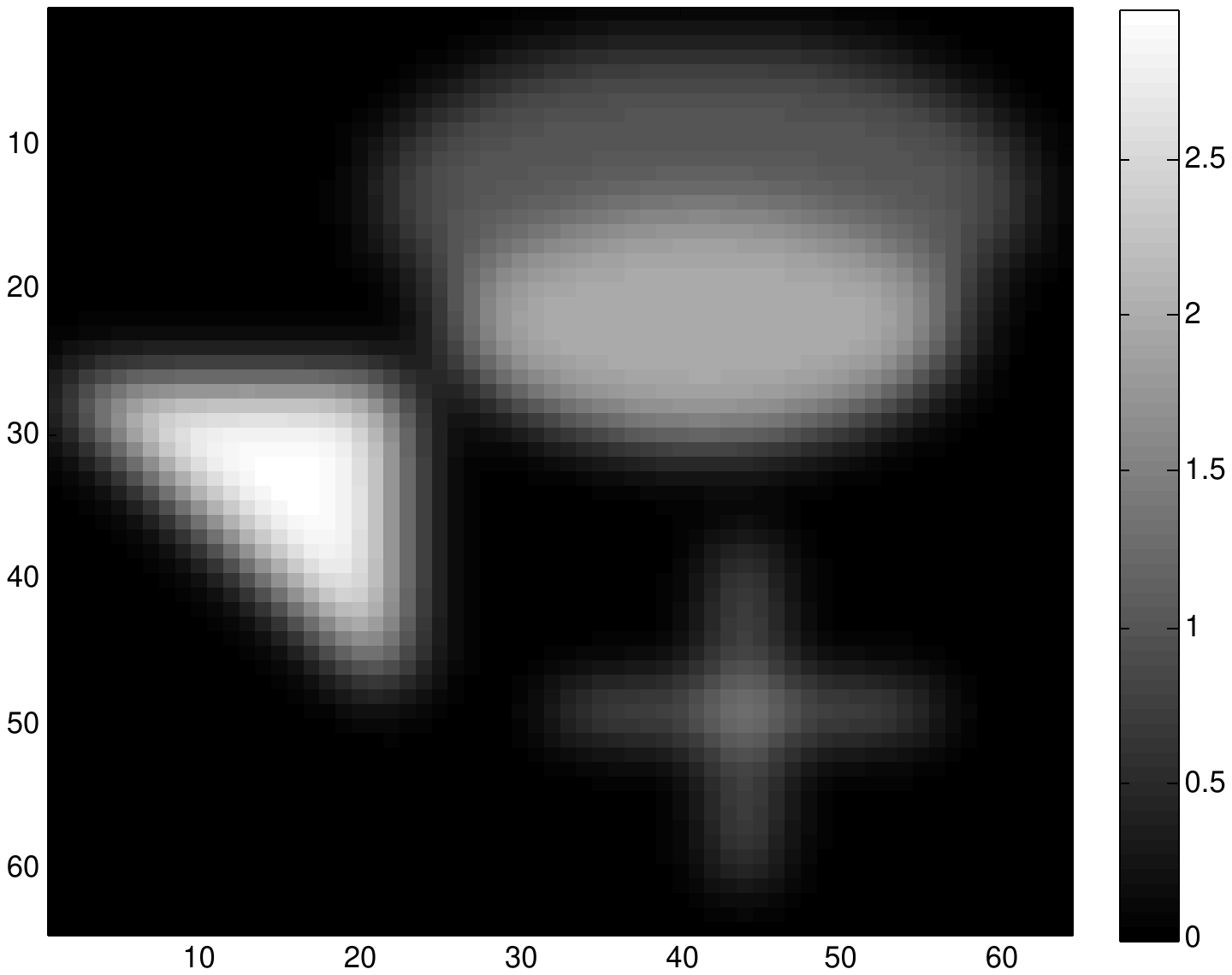} &
\includegraphics[width=6cm]{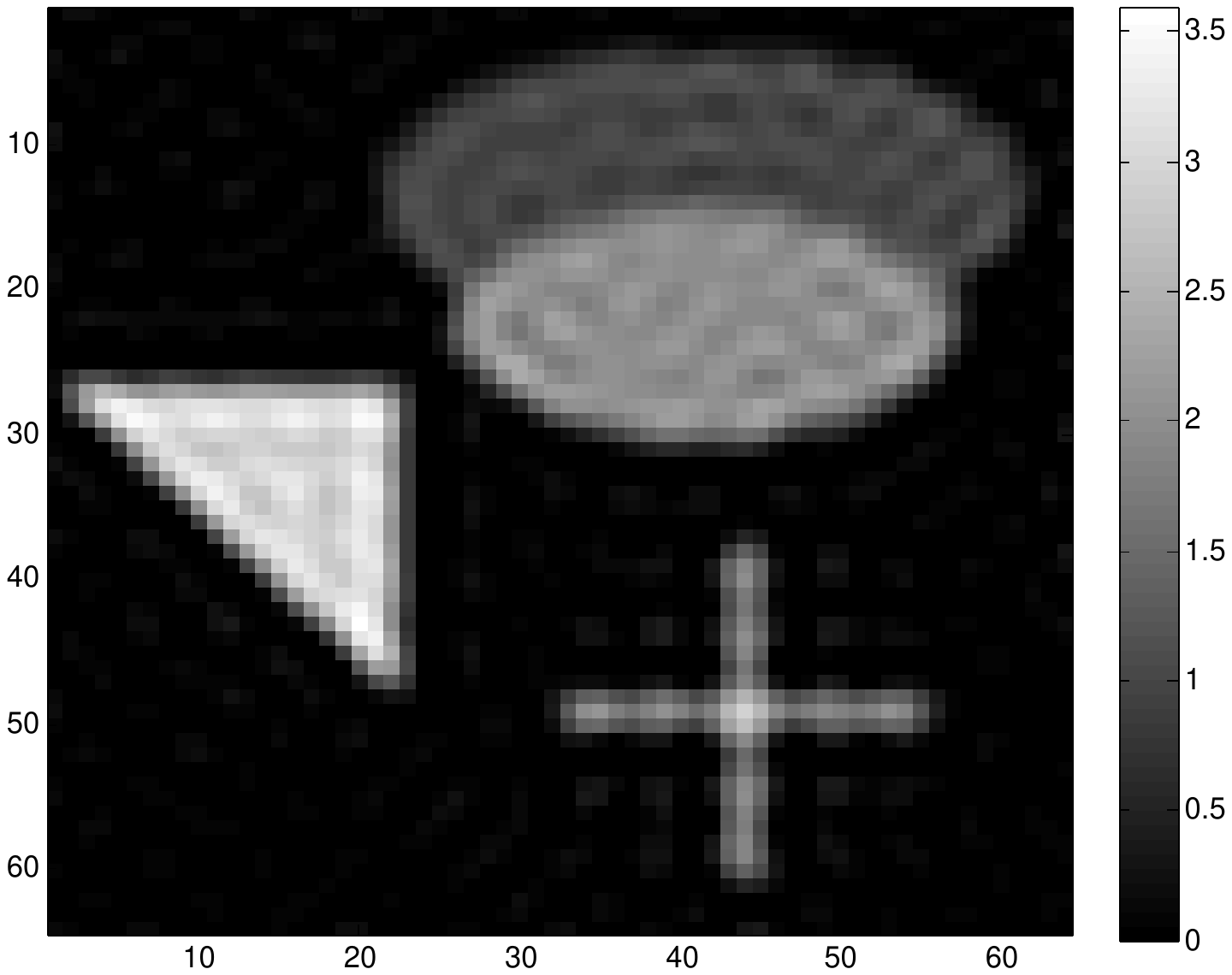} \\
 (a) & (b)
\end{tabular}
\caption{Example 1. Noise level  $0.1\%$. Blurred and noisy image (a) and image restored with the proposed preconditioner (b).}
\label{fig21}
\end{figure}

\begin{figure}[ht]
\centering
\begin{tabular}{cc}
\includegraphics[width=6cm]{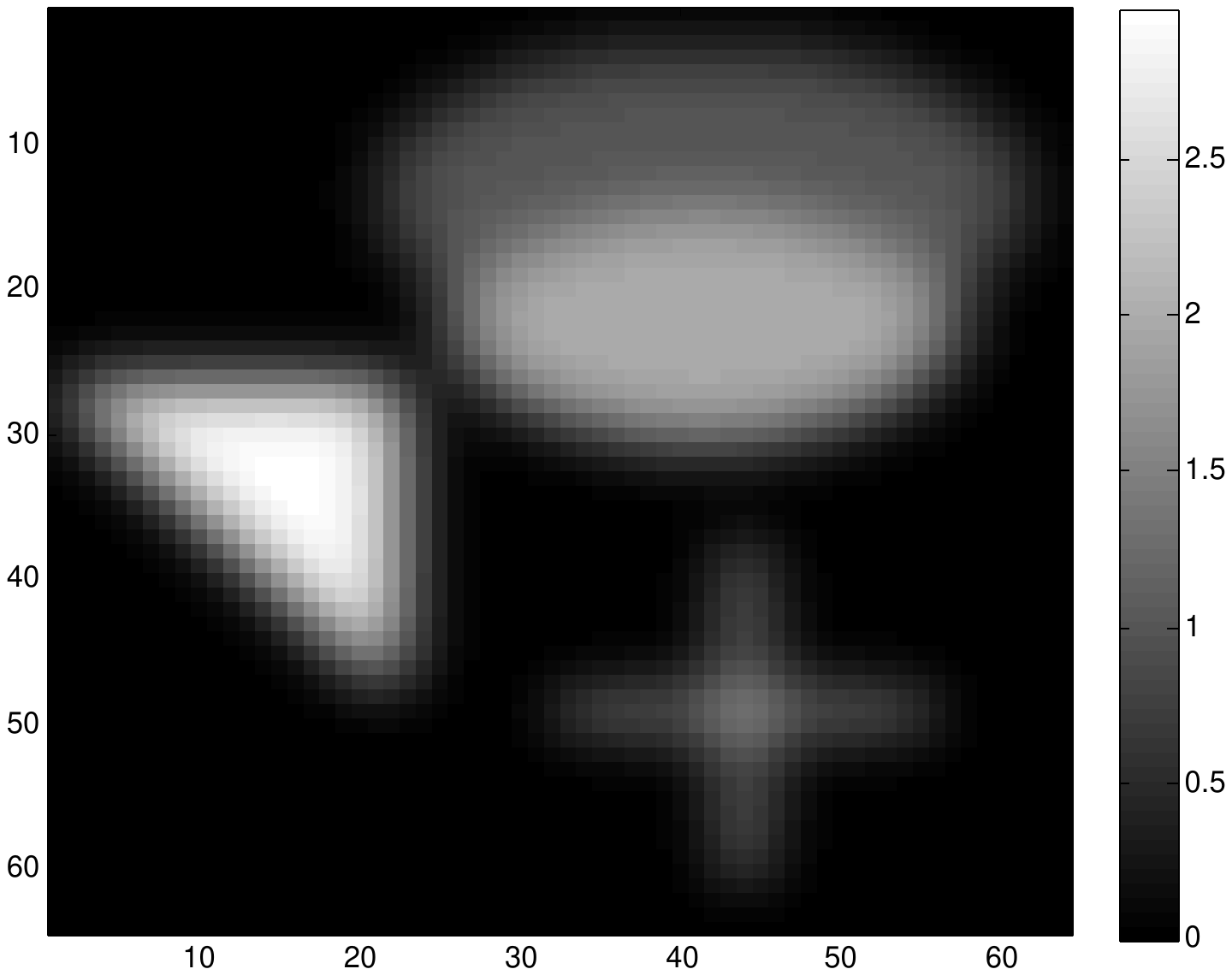} &
\includegraphics[width=6cm]{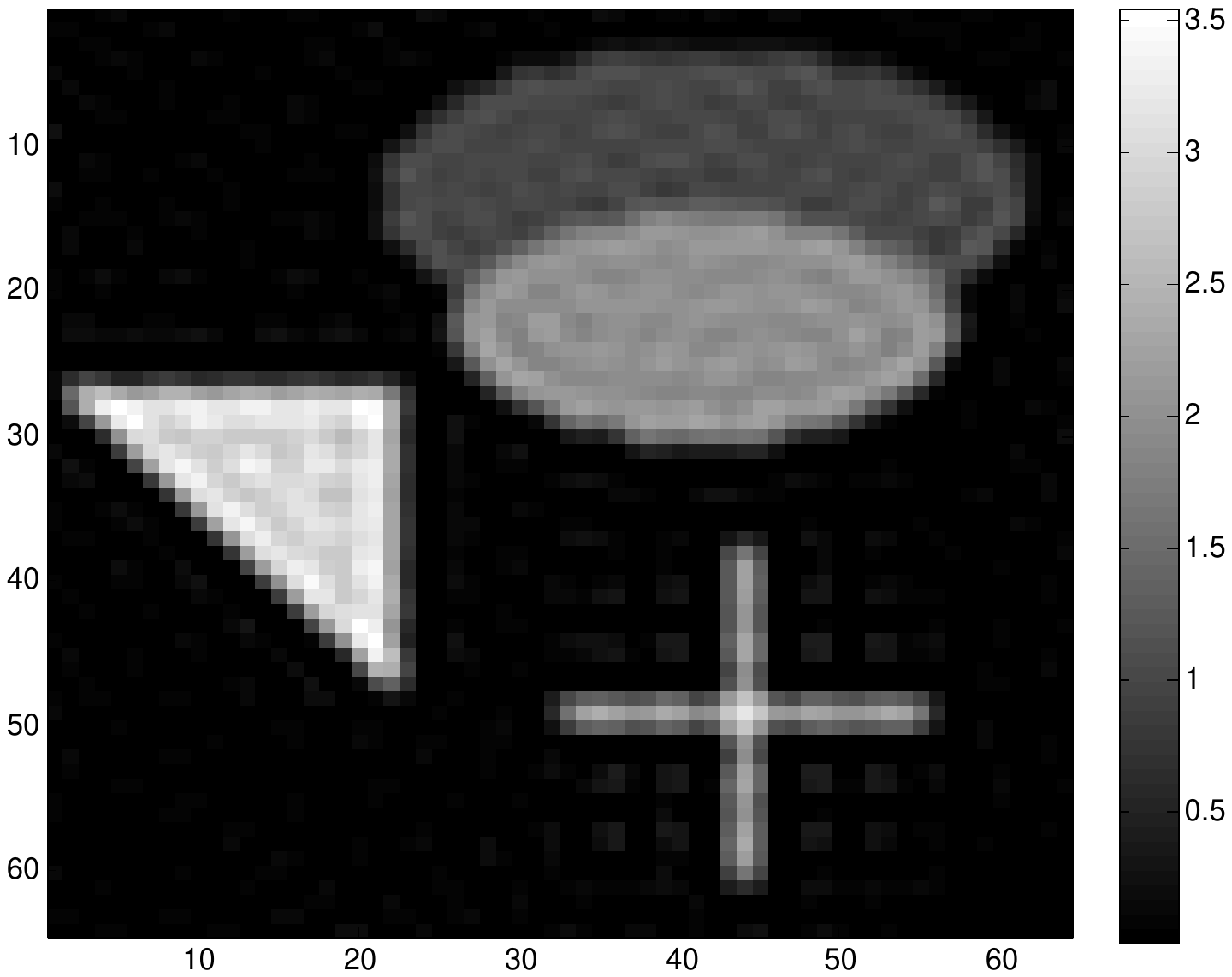} \\
 (a) & (b)
\end{tabular}
\caption{Example 1. Noise level  $0.01\%$. Blurred and noisy image (a) and  image restored with the proposed preconditioner (b).}
\label{fig23}
\end{figure}

We add to the blurred image determined by {\sf blur.m} a noise vector $\be$ with normally 
distributed random entries with mean zero. The vector is normalized to correspond to a 
specified noise level. The pixels of the noise- and blur-contaminated image are ordered 
column-wise and stored in the right-hand side vector $\bb\in\R^{64^2}$. First consider 
$0.1\%$ noise.  Then the parameter $p$ for the proposed
BCCB preconditioner has the value $14$; it is defined using (\ref{choice2_p}). The discrepancy principle prescribes $k=18$ iterations. This yields a restored image with relative error
$0.3404$. Figure \ref{fig21} displays both the available blur- and noise-contaminated image
and the computed restoration. When no preconditioner is used, the discrepancy principle 
terminates the iterations after $k=33$ steps. The restoration so obtained has relative error
$0.3361$. It cannot be distinguished visually from the restoration determined by preconditioned
iterations. We therefore do not display the former. We conclude that preconditioning reduces the
 number of iterations and therefore the computational effort by more than a half and gives a
restoration of about the same quality as unpreconditioned iterations. Table \ref{tab02} displays
the $p$-values used and the number of iterations required for $0.1\%$, $0.05\%$, and 
$0.01\%$ noise in $\bb$.
The noise- and blur-contaminated image together with the restoration determined by preconditioned
iterations for the smallest noise level are displayed in Figure \ref{fig23}.

\begin{figure}[ht]
\centering
\begin{tabular}{cc}
\includegraphics[width=6cm]{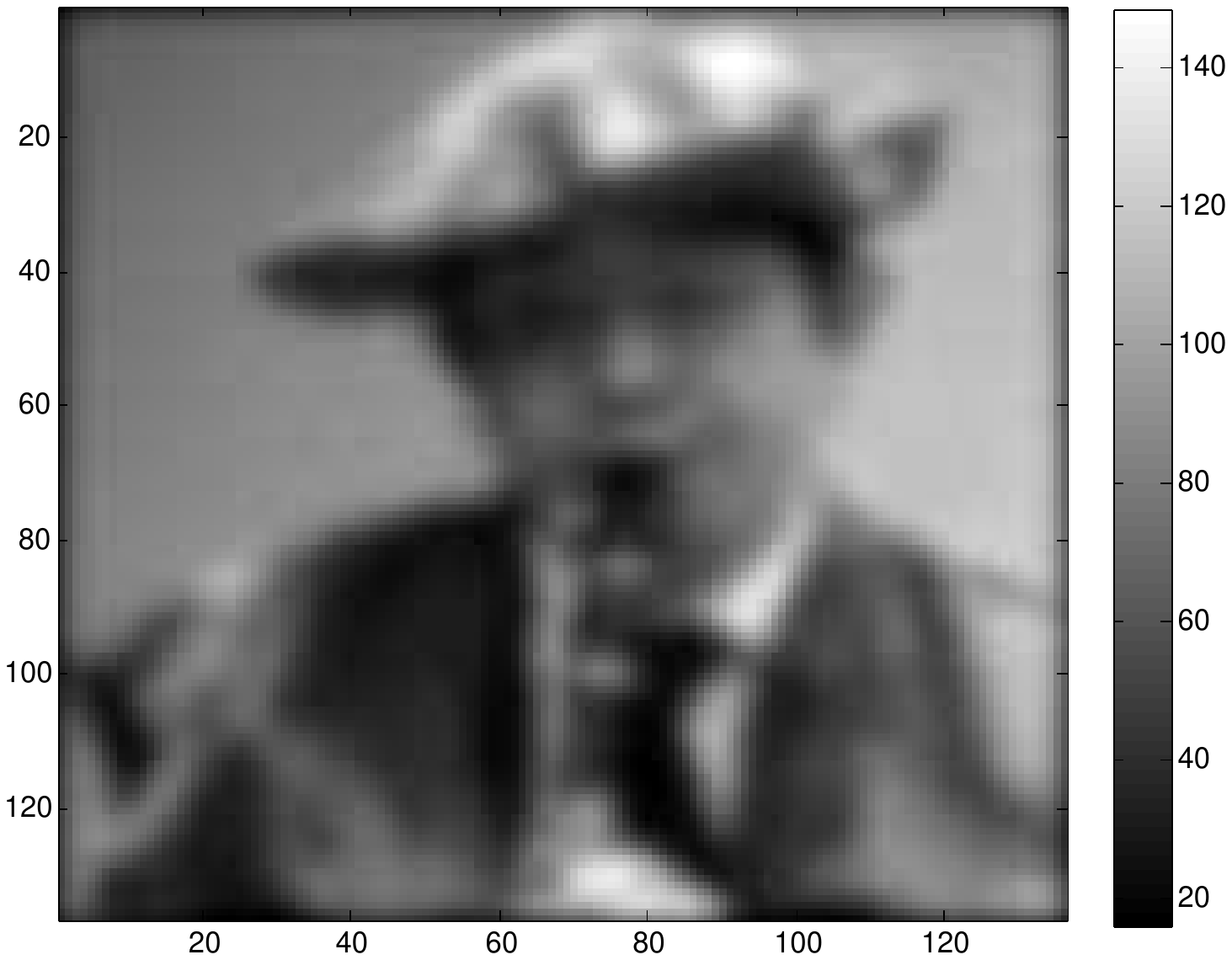} &
\includegraphics[width=6cm]{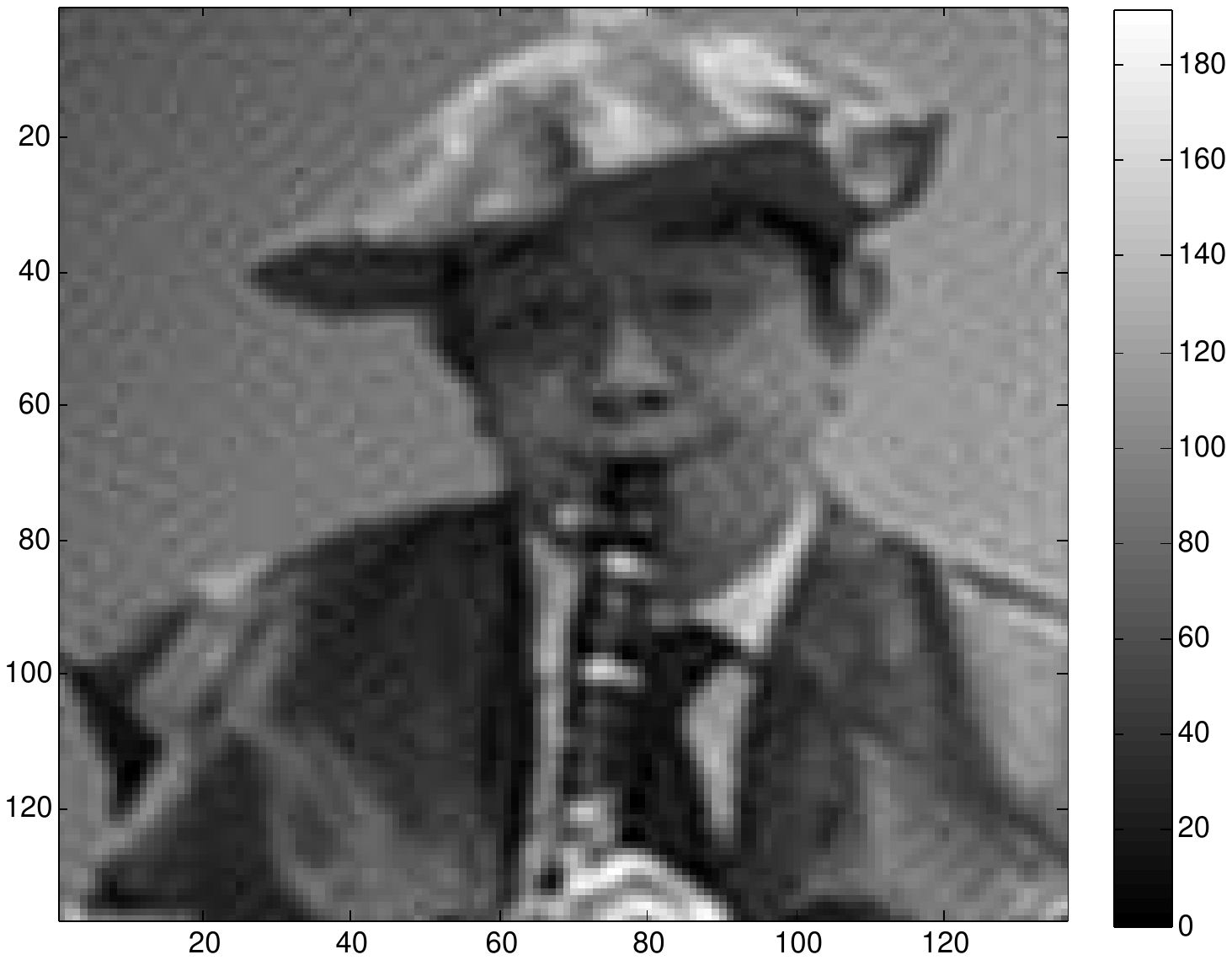} \\
 (a) & (b)
\end{tabular}
\caption{Example 2. Noise level $0.1\%$. Blurred and noisy image (a) and  image restored with the proposed preconditioner (b).}
\label{fig4}
\end{figure}

\begin{figure}[ht]
\centering
\begin{tabular}{cc}
\includegraphics[width=6cm]{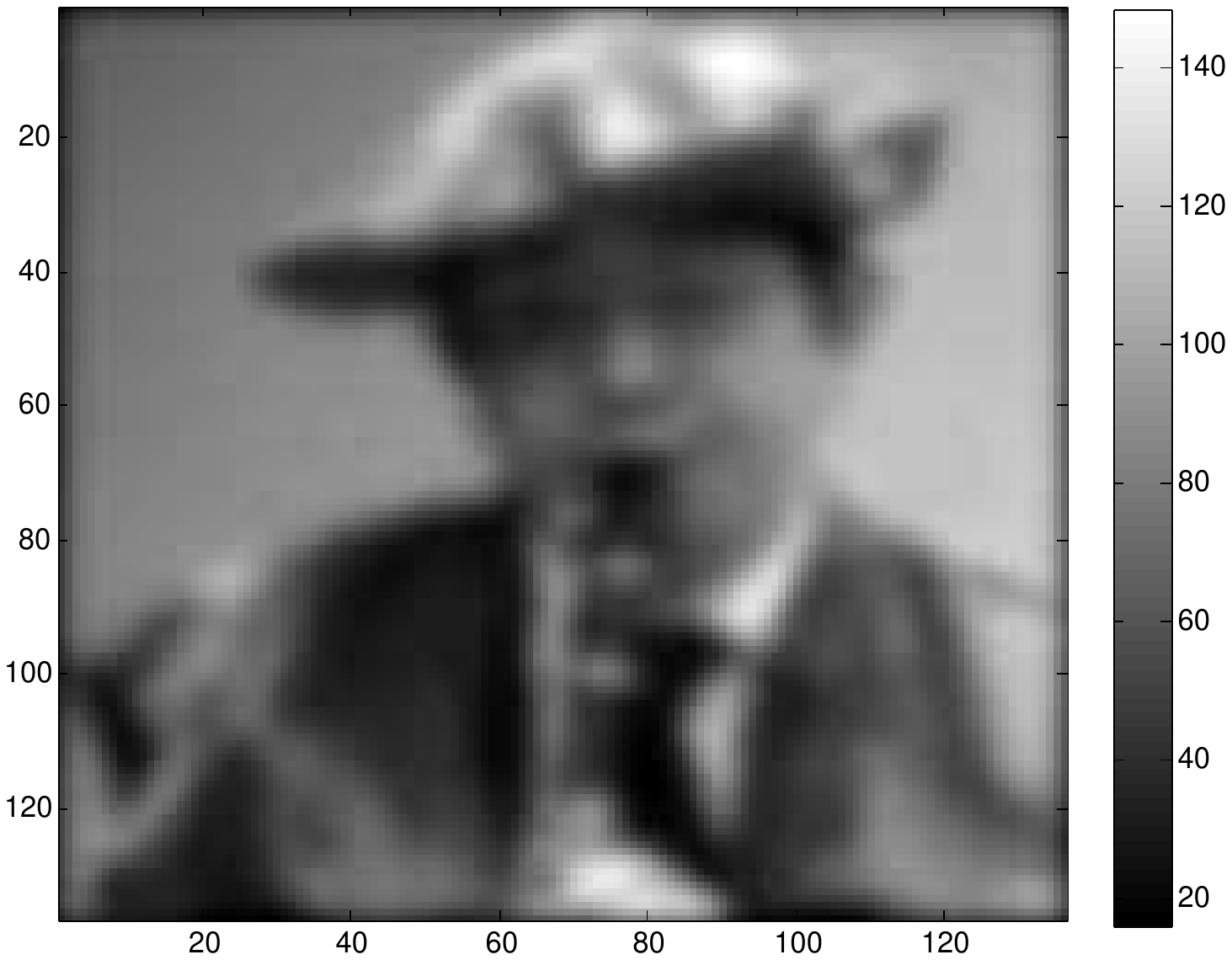} &
\includegraphics[width=6cm]{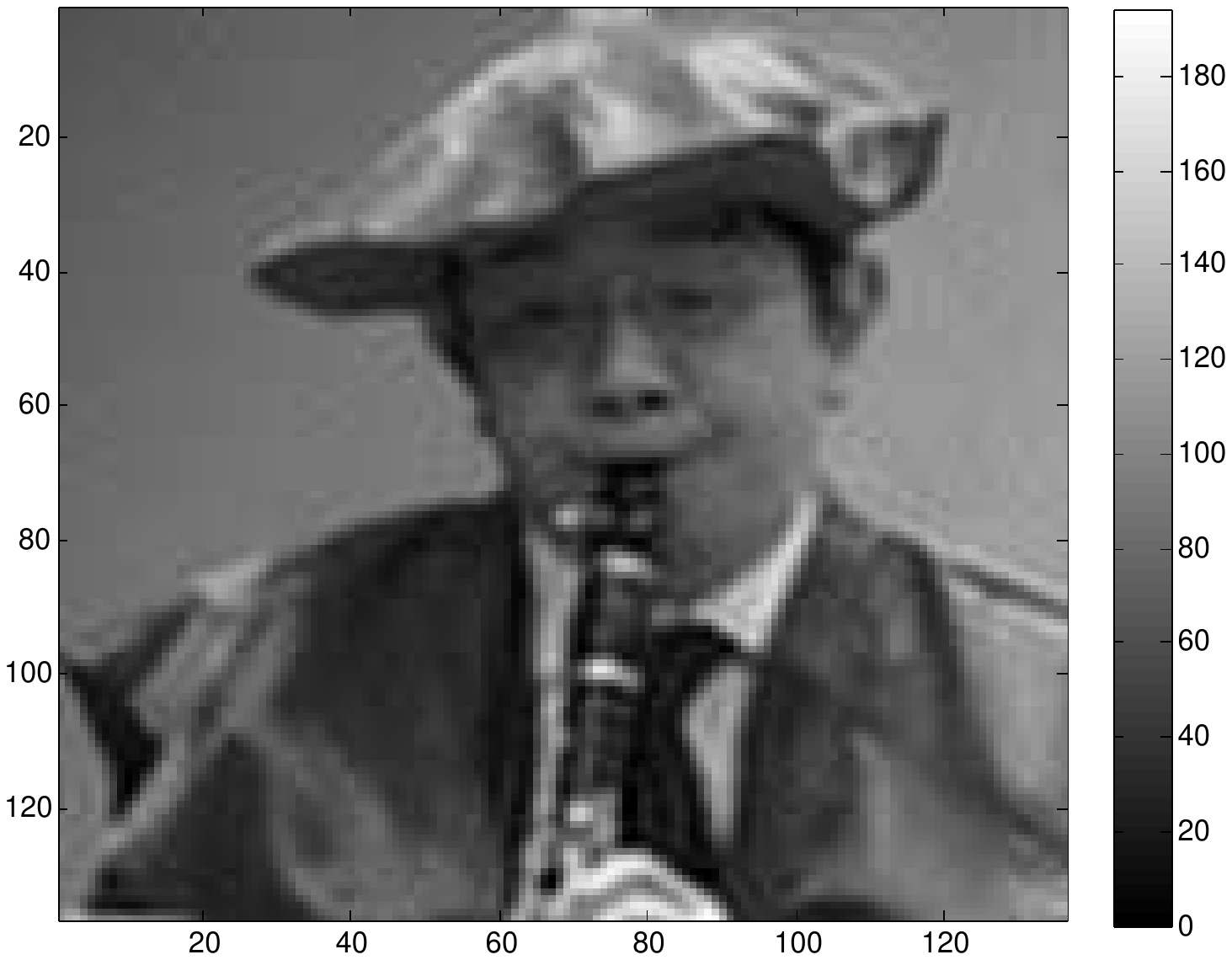} \\
 (a) & (b) 
\end{tabular}
\caption{Example 2. Noise level $0.01\%$. Blurred and noisy image (a) and  image restored with the proposed preconditioner (b).}
\label{fig6}
\end{figure}

{\bf Example 2.} We use the same blur and relative noise as in the previous example, but 
now use the test image ``Ken''. For this image, $n_1=n_2=136$. The BTTB matrix $T$ was
generated by the MATLAB function {\sf blur.m} from \cite{Tools} with the same parameter 
values as in the previous example. We add $0.1\%$, $0.05\%$, and $0.01\%$ white Gaussian 
noise to the blurred image to obtain a blur- and noise-contaminated image, which is stored 
in the right-hand side vector $\bb\in\R^{136^2}$. Table \ref{tab03} displays the number of
iterations required with and without preconditioner to satisfy the discrepancy principle
and the $p$-values that define the preconditioners for different noise levels. 
Figure \ref{fig4} displays the contaminated and restored images
for the noise level $0.1\%$ and Figure \ref{fig6} shows the contaminated and restored images
for the noise level $0.01\%$. Similarly as for Example 1, preconditioned and unpreconditioned
iterations give restorations of essentially the same quality. We therefore only show the
restoration determined by preconditioned iterations. 
 
\begin{table}[htb!]
\centering
\begin{tabular}{cccc}\hline
\% relative data error & $p$ & steps $k$ & $\|\bx_k-\widehat{\bx}\|/\|\widehat{\bx}\|$ \\
\hline
$0.10$ & $27$ & $16$ & $0.0764$ \\
$0.10$ & $-$ & $21$ & $0.0726$ \\ \hline
$0.05$ & $30$ & $16$ & $0.0691$ \\
$0.05$ & $-$ & $28$ & $0.0673$ \\ \hline
$0.01$ & $35$ & $24$ & $0.0571$ \\
$0.01$ & $-$ & $54$ & $0.0575$ \\ \hline
\end{tabular}
\caption{Example 2: ``Ken''. A hyphen signifies that no preconditioner is used.}
\label{tab03}
\end{table}

Finally, consider the situation when the zero vector is chosen as initial approximate 
solution for the preconditioned iterations instead of the vector \eqref{init_guess}.
The parameter $p$ that defines the 
preconditioner is given by (\ref{choice2_p}) and the iterations are terminated by the 
discrepancy principle. Then for noise level $0.1\%$, the discrepancy principle prescribes
$k=20$ iterations and gives a restoration with relative error $0.0904$. The noise level
$0.05\%$ requires $k=21$ iterations and gives a restoration with relative error $0.0776$,
and the noise level $0.01\%$ demands $k=28$ iterations and gives a restoration with 
relative error $0.0596$. A comparison with Table \ref{tab03} shows that the initialization
\eqref{init_guess} requires fewer iterations and gives restorations of higher quality
than when using the initial vector ${\bx}_0={\bzero}$.

{\bf Example 3.}
Our last example is the problem {\sf gravity} from \cite{Tools}. The linear
system of equations \eqref{linsys} is obtained by discretizing an integral 
equation of the first kind with a space invariant kernel. This yields a
Toeplitz matrix $T\in\R^{256\times 256}$ and right-hand side 
$\widehat{\bb}\in\R^{256}$ to which we add an error vector 
${\be}\in\R^{256}$ to obtain the right-hand side of \eqref{linsys}; see 
\eqref{rhs}. The error vector has normally distributed entries with mean 
zero and is scaled to correspond to the noise levels $0.1\%$, $0.05\%$, or 
$0.01\%$. The noise level $0.1\%$ gives the parameter $p=3$ for the circulant
preconditioner and the discrepancy principle is satisfied after $k=8$ 
iterations. We obtain the approximation ${\bx}_8$ of the desired 
solution $\widehat{\bx}$ with relative error $0.0144$. Without a 
preconditioner, the same number of iterations are required to satisfy the 
discrepancy principle and the approximate solution obtained has a larger
relative error, namely $0.0160$. Table \ref{tab01} summarizes the 
results for all noise levels considered. In this example, the preconditioner
does not reduce the number of iterations required to satisfy the discrepancy
principle, but improves  the quality of the computed solution.

In the computations reported in Table \ref{tab01}, we used the initial
iterate \eqref{init_guess}. If instead the initial iterate $\bx_0={\bzero}$
is used for the preconditioned iterations, then the discrepancy
principle is for the noise level $0.1\%$ satisfied after $k=9$ iterations 
and gives an approximate solution with relative error $0.0316$. When
the noise level is reduced to $0.05\%$, the discrepancy principle 
prescribes that $k=10$ iterations be carried out, resulting in an 
approximate solution with relative error $0.0222$ and, finally, for 
noise level $0.01\%$, 
$k=11$ iterations are needed to satisfy the discrepancy 
principle and the computed approximate solution has relative error 
$0.0160$. Thus, for all noise levels the initial iterate 
$\bx_0={\bzero}$ requires more iterations and gives approximate
solutions of inferior quality than the initial iterate \eqref{init_guess}.

\begin{table}[htb!]
\centering
\begin{tabular}{cccc}\hline
\% relative data error & $p$ & steps $k$ & $\|\bx_k-\widehat{\bx}\|/\|\widehat{\bx}\|$ \\
\hline
$0.10$ & $3$ & $8$ & $0.0144$ \\
$0.10$ & $-$ & $8$ & $0.0160$ \\ \hline
$0.05$ & $3$ & $9$ & $0.0105$ \\
$0.05$ & $-$ & $9$ & $0.0119$ \\ \hline
$0.01$ & $3$ & $10$ & $0.0077$ \\
$0.01$ & $-$ & $10$ & $0.0078$ \\ \hline
\end{tabular}
\caption{Example 3: gravity. A hyphen signifies that no preconditioner is used.}
\label{tab01}
\end{table}

\section{Conclusion and extension}\label{sec4} 
This paper presents a novel method to determine BCCB preconditioners to be used for
iterative solution of discretized linear ill-posed problem with a BTTB matrix.
The computed examples show that the number of iterations is reduced by roughly a half
when using the proposed preconditioner, whereas the quality of the computed solution
is about the same as without preconditioning.

Also, we would like to mention that instead of using the circulant preconditioners
described, one may use the generalized optimal circulant preconditioners described
in \cite{NR} at the same computational cost. This may be attractive for certain
Toeplitz and BTTB matrices.


\end{document}